\documentclass{amsart}

\newtheorem{theorem}{Theorem}
\newcommand{\bt}{\begin{theorem}}
\newcommand{\et}{\end{theorem}}
\newtheorem{lemma}{Lemma}
\newcommand{\bl}{\begin{lemma}}
\newcommand{\el}{\end{lemma}}
\newtheorem{corollary}{Corollary}
\newcommand{\bc}{\begin{corollary}}
\newcommand{\ec}{\end{corollary}}
\newtheorem{problem}{Problem}
\newcommand{\bprob}{\begin{problem}}
\newcommand{\eprob}{\end{problem}}
\newtheorem*{conjectureNN}{Conjecture}
\newcommand{\bconjNN}{\begin{conjectureNN}}
\newcommand{\econjNN}{\end{conjectureNN}}
\newtheorem*{theoremNN}{Theorem}
\newcommand{\btNN}{\begin{theoremNN}}
\newcommand{\etNN}{\end{theoremNN}}
\newcommand{\beq}{\begin{equation}}
\newcommand{\eeq}{\end{equation}}
\newcommand{\benum}{\begin{enumerate}}
\newcommand{\eenum}{\end{enumerate}}
\newcommand{\N}{\ensuremath{ \mathbf N }}
\newcommand{\Z}{\ensuremath{\mathbf Z}}

\newcommand{\R}{\ensuremath{\mathbf R}}

\newcommand{\bmat}{\left(\begin{matrix}}
\newcommand{\emat}{\end{matrix}\right)}
\DeclareMathOperator{\qand}{\quad\text{and}\quad}
\DeclareMathOperator{\qqand}{\qquad\text{and}\qquad}

\thanks{Supported in part by a grant from the PSC-CUNY Research Award Program.}

\begin{document}

\title{Cantor  polynomials and the Fueter-P\' olya theorem}
\author{Melvyn B. Nathanson}\address{Department of Mathematics\\Lehman College (CUNY)\\Bronx, NY 10468}\email{melvyn.nathanson@lehman.cuny.edu}

\maketitle

\begin{abstract}
A packing polynomial is a polynomial that maps the set $\mathbf{N}_0^2$ of lattice points 
with nonnegative coordinates bijectively onto $\mathbf{N}_0$.  
Cantor constructed two quadratic packing polynomials, and Fueter and P\' olya 
proved analytically 
that the Cantor polynomials are the only quadratic packing polynomials.   
The purpose of this paper is to present a beautiful elementary 
proof of Vsemirnov of the Fueter-P\' olya theorem.  
It is a century-old conjecture that the Cantor polynomials are the \emph{only} packing polynomials 
on $\mathbf{N}_0^2$.  
\end{abstract}

\section{Storing and packing lattice points}
Let $\Z^m$ be the group of lattice points in $\R^m$.  
Let $\N_0 = \{0,1,2,\ldots\}$ be the set of nonnegative integers, and let 
\[
\N_0^m =  \{ (x_1,\ldots, x_m) \in  \Z^m : x_i \in \N_0 \text{ for } i=1,\ldots, m \}
\]
be the additive semigroup of $m$-dimensional nonnegative lattice points.

In computer science there is the problem of saving and retrieving  matrices 
and other multi-dimensional data in linear memory.  
We can think of computer memory as an infinite sequence of ``storage boxes'' 
numbered $0,1,2,3,\ldots$.  
Suppose that we want to store the matrix $(a_{i,j})$ in memory.  
A simple procedure would be to choose an easily computed one-to-one function 
$F: \N_0^2 \rightarrow \N_0$, and store the matrix element $a_{i,j}$ in the memory box 
with the number $F(i,j)$.    
More generally, for any subset $S$ of $\Z^m$, a one-to-one function 
from $S$ into $\N_0$ is called a \emph{storing  function} on $S$. 
A function that maps $S$ bijectively onto $\N_0$ is called a \emph{packing function} on $S$.  
Computer scientists are interested in finding polynomial 
or other ``elementary'' storing and packing functions on $S$.

For example, let $m=2$ and consider the function $C_1: \N_0^2 \rightarrow \N_0$ 
that enumerates, from lower right to upper left,  
the nonnegative lattice points on the consecutive parallel line 
segments $x+y = k$ for $k=0,1,2,\ldots$.  
Thus, $C_1(0,0) = 0$, $C_1(1,0) = 1$, $C_1(0,1) = 2$, $C_1(2,0) = 3,$ $C_1(1,1) = 4$, 
$C_1(0,2)=5$, $C_1(3,0)=6$, \ldots.
For $k \geq 0$,  the number of nonnegative lattice points $(x,y)$ with $x+y < k$ is $1 + 2 + \cdots + k = k(k+1)/2$.  
It follows that $C_1(k,0) = k(k+1)/2$.  If $(x,y) \in \N_0^2$ and $x+y=k$, then 
\begin{align*}
C_1(x,y) 
& = C_1(k,0) + y = C_1(x+y,0) + y  \\
& = \frac{(x+y)(x+y+1)}{2} + y \\
& = \frac{1}{2} \left( x+y \right)^2 + \frac{1}{2} \left( x + 3y \right).
\end{align*}
The function $C_2: \N_0^2 \rightarrow \N_0$ 
that enumerates, from upper left to lower right,  
the nonnegative lattice points on the consecutive parallel line 
segments $x+y = k$ is 
\[
C_2(x,y) 
= \frac{1}{2} \left( x+y \right)^2 + \frac{1}{2} \left( 3x + y \right).
= C_1(y,x).
\]
The functions $C_1(x,y)$ and $C_2(x,y)$  are called 
the \emph{Cantor packing polynomials}.

The classical applications of polynomial packing functions are due to Cauchy and Cantor.
Cantor~\cite[page 494]{cant95},~\cite[page 107]{cant52} 
used  the quadratic polynomials $C_1(x,y)$ and $C_2(x,y)$ 
to prove that the set $\N_0^2$ is countable.  
Earlier, to reduce double sums to single sums, 
Cauchy~\cite[Part 1, Chapter 6, Theorem 6]{cauc21} 
used the equivalent identity 
\[
\sum_{ i=0 }^{\infty} \sum_{ j=0 }^{\infty} a_i b_j =   \sum_{k=0}^{\infty} c_k
\]
where
\[
c_k = \sum_{ \substack{    i,j=0 \\i+j=k}    }^k a_{k-i}b_i 
= a_kb_0 + a_{k-1}b_1 +  a_{k-2}b_2 +\cdots + a_0b_k.  
\]

The following simple argument implies that a  packing function cannot be linear.

\bl               \label{Cantor-FP:lemma:NoLinear}
Let $\ell, m \in \N_0$ with $m \geq 2$.    
For $x = (x_1,\ldots, x_m) \in \N_0^m$, let $\min(x)= \min( x_i:i=1,\ldots, m)$, 
and let $S$ be the set of all lattice points 
$x \in \N_0^m$ such that $ \min(x) \geq \ell$.  
No linear polynomial is a storing function on $S$.
\el

\begin{proof}
For $x = (x_1,\ldots, x_m) \in \R^m$, consider the linear polynomial 
\[
F(x) = F(x_1,x_2,\ldots, x_m) = a_1x_1 + a_2x_2  + \cdots + a_mx_m + c.
\]
Let $A = \max(|a_i| : i=1,\ldots, m)$.  
Recall the Kronecker delta 
\[
\delta_{i,j}
= \begin{cases}
1 & \text{if $i=j$} \\
0 & \text{if $i \neq j$.}
\end{cases}
\]
For $i=1,\ldots, m$, let $e_i = (\delta_{i,1}, \delta_{i,2},\ldots, \delta_{i,m}) \in \N_0^m$. 

Suppose that  $F$ is a storing function on $S$, that is, $F(S) \subseteq \N_0$ and 
$F$ is one-to-one on $S$.
If $x \in S$, then $x+e_i\in S$ and  
\[
F(x+e_i) - F(x) = a_i  \in \Z.
\]
It follows that $\sum_{i=1}^m a_ix_i \in \Z$, 
and so $c = F(x) - \sum_{i=1}^m a_ix_i \in \Z$.

Suppose that $F: S \rightarrow \N_0$ is a sorting function.  
Let $i,j \in \{1,2,\ldots, m\}$ with $i \neq j$.  
If $x \in S$ and 
$\min( x ) \geq \ell + A$, 
then $x+ a_i e_j - a_j e_i \in S$ and 
\[
F(x+a_ie_j-a_je_i) = F(x) + a_ia_j-a_ja_i = F(x).
\]
Because $F$ is one-to-one on $S$, it follows that $x = x+a_ie_j-a_je_i$ 
and so $a_ie_j = a_je_i$.  
This implies that $a_i = a_j = 0$ and so $a_i = 0$ for all $i \in \{1,\ldots, m\}$.
Thus, the storing function $F(x) = c$ is constant, which is absurd.  
This means that no linear polynomial can be a storing function on $S$.   
\end{proof}

\bconjNN
The Cantor packing polynomials $C_1(x,y)$ and $C_2(x,y)$ are the \emph{only} polynomial 
bijections from $\N_0^2$ to $\N_0$.
\econjNN

This conjecture is nearly 100 years old.   
In 1923, Fueter and P\' olya~\cite{fuet-poly23} obtained the first result 
about the uniqueness of packing polynomials in two variables.  
Using methods from analytic number theory, 
they proved that the  Cantor polynomials  are the unique \emph{quadratic} 
polynomial packing functions from $\N_0^2$ to $\N_0$.  
In 1978, Lew and Rosenberg~\cite{lew-rose78b,lew-rose78a} proved that  
no polynomial packing function on $\N_0^2$ has degree three or four.   
It is not known if there exists a packing polynomial on $\N_0^2$ of degree greater than four.  

There exist many packing polynomials for lattice points of dimension $m \geq 3$.  
The simplest are constructed by composing Cantor polynomials.  
For example, the function $(x,y,z) \mapsto (C_1(x,y),z)$ is a bijection 
from $\N_0^3$ to $\N_0^2$, and so 
$(x,y,z) \mapsto C_1(C_1(x,y),z)$ is a packing polynomial on $\N_0^3$.  
There also exist packing polynomials on $\N_0^m$ that are not compositions of  
packing polynomials in lower dimensions.  
Much work has been done on this problem, 
e.g.~\cite{fett-arre-mora05,lew81,mora97,mora16,mora-arre99,mora-lew96,mora-sanc00,rose71,rose78,sanc95}.


In 2001, Vsemirnov~\cite{vsem01,vsem02} gave a beautiful proof 
of the Fueter-P\' olya theorem that uses only Gauss's law of quadratic reciprocity 
and Dirichlet's theorem on primes in arithmetic progressions.   
The purpose of this paper is to present Vsemirnov's  proof.

\section{Elementary proof of  the Fueter-P\' olya theorem}

We begin with an exercise in elementary number theory.  

\bl     \label{Cantor-FP:lemma:nonresidue}
If $D$ and $\ell$ are nonzero integers and $D$ is not a square, 
then there exists a prime $p$ such that $D$ is a 
quadratic non-residue modulo $p$ and $p$ does not divide $\ell.$  
\el

\begin{proof}
We can write 
\[
D = (-1)^{\alpha} 2^{\beta} m^2 \prod_{i=1}^k q_i 
\]
where $\alpha, \beta \in \{0,1\}$, $m \in \N$, and $q_1,\ldots, q_k$ are distinct odd primes.

Suppose that $k=0$. 
If $\beta = 0$, then $\alpha = 1$ and $D = -m^2$.  
Choosing a prime $p \equiv 3 \pmod{4}$, we obtain 
\[
\left(\frac{D}{p}\right) =  \left(\frac{-1}{p}\right) = -1.  
\]
If $\beta = 1$, then $D = \pm 2 m^2$.
Choosing a prime $p \equiv 5 \pmod{8}$, we obtain 
\[
\left(\frac{D}{p}\right) =  \left(\frac{ \pm 2}{p}\right) = -1.  
\]
By Dirichlet's theorem on primes in arithmetic progressions, 
each of these congruence classes contains infinitely many primes.  
Choosing $p > \ell$, we obtain a prime $p$ does not divide $\ell$.

Suppose that $k \geq 1$.  
For every prime number such that $p \equiv 1 \pmod{8}$, we have 
\[
\left(\frac{-1}{p}\right) = \left(\frac{2}{p}\right) = 1.
\]
Applying the multiplicativity of the Legendre symbol and quadratic reciprocity, we obtain 
\[
\left(\frac{D}{p}\right) = \prod_{i=1}^k  \left(\frac{q_i}{p}\right) 
= \prod_{i=1}^k  \left(\frac{p}{q_i}\right).  
\]
Choose integers $r_1,\ldots, r_k$ such that 
\[
\left(\frac{r_1}{q_1}\right) = -1
\]
and 
\[
\left(\frac{r_i}{q_i}\right) = 1 \qquad \text{for $i=2,\ldots, k$.}
\]
By the Chinese remainder theorem, there is an integer $s$ such that 
\[
s \equiv 1 \pmod{8} 
\]
and
\[
s \equiv r_i \pmod{q_i}  \qquad \text{for $i=1,2,\ldots, k$.}
\]
Moreover,
\[
\left(s, 8 \prod_{i=1}^k q_i \right) = 1.
\]
If $p$ is a prime number such that 
\beq         \label{Cantor-FP:congruence2}
p \equiv s \pmod{8 \prod_{i=1}^k q_i}
\eeq
then
\[
\left(\frac{D}{p}\right) =  \left(\frac{q_1}{p}\right) = -1.  
\]
By Dirichlet's theorem on primes in arithmetic progressions, 
there are infinitely many primes $p$ that satisfy the congruence~\eqref{Cantor-FP:congruence2}.
Choosing a prime $p$ in this arithmetic progression such that $p > \ell$ 
completes the proof.  
\end{proof}

We can now begin the proof of Fueter-P\' olya theorem.  There are four lemmas.

\bl   \label{Cantor-FP:lemma:Quadratic}
If $F(x,y)$ is a quadratic packing polynomial, then there exist nonnegative 
integers $a,c,f$ and integers $b,d,e$  such that 
\[
a\equiv d \pmod{2}   
\]
\[
c\equiv e \pmod{2}
\]
and 
\beq     \label{Cantor-FP:standardF}
F(x,y) =  \frac{1}{2}\left( ax^2 +2bxy + cy^2 \right) +    \frac{1}{2}\left(dx+ey\right) + f.
\eeq   
Moreover, if $a=c=0$, then $b \geq 1$.  
\el

\begin{proof}
Every quadratic polynomial $F(x,y)$ with complex coefficients can be written in the 
form~\eqref{Cantor-FP:standardF}.  
Because $F$ is a function from $\N_0^2$ to $\N_0$, we have 
\[
a = F(2,0)-2F(1,0)+F(0,0)
\in \Z
\]
\[
c = F(0,2)-2F(0,1)+F(0,0) \in \Z
\]
and
\[
f =  F(0,0) \in \N_0.
\]
For all $x,y \in \N_0$, the inequalities 
\[
F(x,0) =  \frac{a}{2} x^2  + \frac{d}{2}x  + f \geq 0
\]
and
\[
F(0,y)  = \frac{c}{2} y^2 +  \frac{e}{2} y +  f \geq 0
\]
imply that $a$ and $c$ are nonnegative integers.  
The identities 
 \[
F(1,0) - F(0,0) = \frac{a+d}{2}  \in \Z
\]
 \[
F(0,1) - F(0,0) = \frac{c+e}{2}  \in \Z
\]
imply that $d$ and $e$ are integers such that 
\[
a \equiv d \pmod{2} \qqand c \equiv e \pmod{2}.
\]
It follows that 
\[
 F(1,1) = b+  \frac{a+d}{2} + \frac{c+e}{2} + f \in \N_0
\]
and so $b$ is an integer.  

If $a = c = 0$, then $b \neq 0$ because $F$ is quadratic.  
For all $x \in \N_0$,  
\[
F(x,x) = bx^2 +  \frac{ d+e}{2} x + f \geq 0
\]
and so  $b \geq 1$.
This completes the proof.
\end{proof}

\bl             \label{Cantor-FP:lemma:QuadraticForm}
If $F(x,y)$ is a quadratic packing polynomial of the form~\eqref{Cantor-FP:standardF}, 
then the quadratic form 
\[
Q(x,y) = \frac{1}{2} \left( ax^2 + 2bxy + cy^2 \right) 
\]
is positive-definite on $\N_0^2$.   Moreover, $a \geq 1$ and $c \geq 1$.
\el

\begin{proof}
The quadratic form  $Q(x,y)$ is  nonzero because the polynomial $F(x,y)$ is quadratic.
Defining the linear form 
\[
L(x,y) =  \frac{1}{2}(dx+ey)
\]
we can write
\[
F(x,y) = Q(x,y) + L(x,y) + f
\]
If $a \geq 1$ and $r > |d|/2$, then $Q(r,0) > |L(r,0)|$.  
If $c \geq 1$ and $s > |e|/2$, then $Q(0,s) > |L(0,s)|$.  
If $a=c=0$ and $r > (|d|+|e|)/2$, then $b \geq 1$ by Lemma~\ref{Cantor-FP:lemma:Quadratic},
and $Q(r,r) > |L(r,r)|$.   
Thus, there exists $(r,s) \in \N_0^2 \setminus \{ (0,0)\}$ such that 
\beq        \label{Cantor-FP:Qrs}
Q(r,s) > |L(r,s)|.
\eeq

It is easy to show that $Q$ is nonnegative-definite.  
For all $(x,y) \in \N_0^2$ and  $t \in \N_0$, 
\[
F(xt,yt) = Q(xt,yt) +L(xt,yt) + f 
= Q(x,y)t^2 + L(x,y)t + f.
\]
If $Q(x,y) < 0$ for some $(x,y) \in \N_0^2$, 
then $F(xt,yt)  < 0$
for all sufficiently large $t$, which is absurd.
Therefore, $Q(x,y) \geq 0$ for all $(x,y) \in \N_0^2$.

Suppose that $Q(u,v) = 0$ for some $(u,v) \in \N_0^2 \setminus \{ (0,0)\}$.  
For all $t \in \N_0$, we have 
\[
F(ut, vt) = Q(u,v)t^2 + L(u,v)t + f = L(u,v)t + f = w t + f
\]
where $w = L(u,v)$ is a positive integer because $F(ut, vt)$ is nonconstant 
and nonnegative for $t \in \N_0$.    
Choosing $(r,s) \in \N_0^2$ that satisfies inequality~\eqref{Cantor-FP:Qrs}, we have  
\[
F(rw, sw) = Q(r,s)w^2 + L(r,s)w + f = w  m + f
\]
where
\[
m =  Q(r,s) w+ L(r,s) \geq Q(r,s) - |L(r,s)| > 0.
\]
Because $F$ is one-to-one on $\N_0^2$ and 
\[
F(um,vm) = w m + f = F(rw, sw)
\]
it follows that 
\[
(um,vm) = (rw, sw)
\]
and so 
\[
0 = Q(u,v) m^2 = Q(um,vm)  = Q(rw, sw) = Q(r,s) w^2 > 0
\]
which is absurd.  Therefore, the quadratic form $Q(x,y)$ is positive-definite on $\N_0^2$.
This implies that 
\[
a = 2Q(1,0) \geq 1 
\]
and
\[
c = 2Q(0,1) \geq 1.
\]
This completes the proof.  
\end{proof}

\bl                \label{Cantor-FP:lemma:abc}
If $F(x,y)$ is a quadratic packing polynomial of the form~\eqref{Cantor-FP:standardF}, 
then $b \leq 1$.  
\el

\begin{proof}
Let $m \geq \max(2, |d|, |e|)$.  If $b \geq 2$, then 
\begin{align*}
F(x,y) 
& =  \frac{1}{2}\left( ax^2 +2bxy + cy^2 \right) +  \frac{1}{2}\left( dx+ey\right) + f \\
& \geq  \frac{1}{2}\left( x^2 +4xy + y^2 \right) - \frac{1}{2}\left( |d|x+ |e|y\right)  \\
& \geq  \frac{1}{2}(x+y)^2 + xy - \frac{m}{2}\left(x+ y\right).
\end{align*}
We obtain the following inequalities:
\beq           \label{Cantor-FP:region-1}
F(x,y) \geq  \frac{1}{2} (x+y) (x+y-m) + xy   
\eeq
and, for every positive integer $k$, 
\beq  \label{Cantor-FP:region-3}
F(x,y) 
=  \frac{1}{2} (x+y)^2 + \frac{k-1}{k} xy  + \frac{1}{2k} (x-km)y + \frac{1}{2k} (y-km)x.
\eeq

There are  three cases.  \\
Case 1.  
If $\max(x,y) \geq 25 m$, then $x+y \geq 25m$.  
Applying inequality~\eqref{Cantor-FP:region-1}, we obtain  
\[
F(x,y) \geq  \frac{1}{2} (x+y) (x+y-m) \geq  \frac{1}{2} (25m) (24m) =  300m^2.
\]
Case 2. 
If $\min(x,y) \geq 10 m$,  then $x+y \geq 20m$.  
Applying inequality~\eqref{Cantor-FP:region-3} with $k=10$,
we obtain  
\[
F(x,y) \geq  \frac{1}{2} (20m)^2+ \frac{9}{10}  (10m)^2 = 290 m^2.
\]
Case 3.  
If $\min(x,y) \geq m$ and $x+y \geq 24m$, then, applying inequality~\eqref{Cantor-FP:region-3} 
with $k=1$, we obtain 
\[
F(x,y) \geq  \frac{1}{2}  (24m)^2  =288m^2.  
\]
It follows that if $(x,y) \in \N_0^2$  and $F(x,y)  < 288 m^2$, then 
\[
\max(x,y) < 25 m
\]
\[
\min(x,y) < 10 m
\]
and
\[
\min(x,y) < m  \quad \text{ or } \quad x+y < 24m,
\]
Equivalently,  the lattice point $(x,y)$ must belong to exactly one 
of the following five sets:
\begin{align*}
Z_1 & = \left\{ (x,y) \in \N_0^2 :   0 \leq x < m \qand 0 \leq y < 25m   \right\}  \\
Z_2 &= \left\{ (x,y) \in \N_0^2 :  m \leq x < 10 m \qand 0 \leq y < 24m - x     \right\}  \\
Z_3 & = \left\{ (x,y) \in \N_0^2 :  10 \leq x < 14 m \qand 0 \leq y < 10m    \right\}   \\
Z_4 & = \left\{ (x,y) \in \N_0^2 :  14 \leq x < 23 m \qand 0 \leq y < 24m - x    \right\}  \\  
Z_5 & = \left\{ (x,y) \in \N_0^2 :  23 m \leq x < 25 m \qand 0 \leq y < m    \right\}.  
\end{align*}
For $i = 1,\ldots, 5$, let $N_i$ denote the number of lattice points in the set $Z_i$.  
We have
\begin{align*}
N_1 & = 25 m^2 \\
N_2 & = \frac{333}{2} m^2 +  \frac{9}{2}m \\
N_3 & = 40 m^2 \\
N_4 & =  \frac{99}{2} m^2 +  \frac{9}{2}m \\
N_5 & = 2 m^2.
\end{align*}
Therefore, the number of nonnegative integers $n < 288 m^2$ 
represented by the polynomial $F(x,y)$ with $(x,y) \in \N_0^2$ is at most 
\[
\sum_{i=1}^5 N_i = 283m^2 + 9m < 288m^2
\]
because $m \geq 2$.  
It follows that $F:\N_0^2 \rightarrow \N_0$ is not surjective, which is absurd.
Therefore, $b \leq 1$.
This completes the proof. 
\end{proof}

\bl                 \label{Cantor-FP:lemma:CompleteSquare}
If $F(x,y)$ is a quadratic packing polynomial of the form~\eqref{Cantor-FP:standardF}, 
and if 
\[
D = b^2-ac 
\]
then 
\[
8 aD F(x,y) =  D u^2  - v^2  + r 
\]
where 
\beq                   \label{Cantor-FP:vse-1}
u = 2ax + 2b y  + d 
\eeq
\beq                   \label{Cantor-FP:vse-2}
v =   2D y + bd - ae 
\eeq
and
\beq                   \label{Cantor-FP:vse-3}
r =    (bd - ae)^2 - D d^2 + 8aDf. 
\eeq
\el

\begin{proof}
This is an exercise in ``completing the square.''
We have 
\begin{align*}
8 aD F(x,y) & - 8 aD f =  4D (a^2 x^2 + 2abxy + acy^2 + adx + aey)   \\
& =  4D \left( (ax+by)^2 -  D y^2 + adx + aey  \right)   \\
& =  4D \left( (ax+by)^2+ d(ax+by)  -  D y^2 -  ( bd - ae) y  \right)   \\
& =  D \left( 4(ax+by)^2+ 4d(ax+by) \right)  - \left( 4D^2 y^2 + 4D  ( bd - ae) y  \right)   \\
& =  D \left( 2ax+2by+ d\right)^2  - \left( 2D y +  bd - ae  \right)^2  +  \left(bd - ae  \right)^2  
- D d^2 \\
& = Du^2 - v^2 +r - 8aDf.
\end{align*}
This completes the proof.  
\end{proof}

We can now complete Vsemirnov's proof of the Fueter-P\' olya theorem.

\btNN
Every quadratic packing polynomial is a Cantor polynomial.
\etNN

\begin{proof}
We shall start by proving that $D = b^2 - ac$ is a square.  
Suppose not.  Then $D  \neq 0$.  
By Lemma~\ref{Cantor-FP:lemma:CompleteSquare}, 
if $x$ and $y$ are integers, then there are integers $u = u(x,y)$, $v = v(y)$, and $r$ such that 
\[
8 aD F(x,y) =    D u^2  - v^2  + r.
\]
Because  $8a\neq 0$ and $D$ is not a square, 
we apply Lemma~\ref{Cantor-FP:lemma:nonresidue} with 
$\ell = 8a$ and obtain a prime number $p$ such that 
$D$ is a quadratic non-residue modulo $p$ 
and $(8 aD,p) = 1$.  
There is an integer $s$ such that $(s,p) = 1$ and 
\[
8 aD s \equiv  r \pmod{p}.
\]
Because $F(x,y)$ is a packing polynomial, 
there exist infinitely many lattice points $(x,y) \in \N_0^2$ such that 
$F(x,y) \equiv s \pmod{p}$ and 
\[
 D u^2  - v^2  + r = 8 aD F(x,y) \equiv 8 aD s \equiv  r \pmod{p}.  
\]
Therefore,
\[
 D u^2  \equiv v^2  \pmod{p}.  
\]
If $u \not\equiv 0 \pmod{p}$, then 
\[
D \equiv \left( v u^{-1} \right)^2 \pmod{p}
\]
which is impossible because $D$ is a quadratic non-residue modulo $p$.  
It follows that  $u \equiv 0 \pmod{p}$,  and so $ v \equiv 0 \pmod{p}$ and 
\[
D u^2  - v^2  \equiv 0 \pmod{p^2}.
\]
It follows that if $(x,y) \in \N_0^2$ and $F(x,y) \equiv s \pmod{p}$, then 
\[
8 aD  F(x,y) = D u^2  - v^2 + r  \equiv r \equiv 8 aD s \pmod{p^2}.
\]
Because $(8aD,p)=1$, we see that the congruence $F(x,y) \equiv s \pmod{p}$ 
implies that $F(x,y) \equiv s \pmod{p^2}$.  
This means that there do not exist integers $x$ and $y$ such that 
$F(x,y) \equiv s+p \pmod{p^2}$, and so the polynomial $F(x,y)$ is not surjective 
from $\N_0^2$ onto $\N_0$, which is absurd.

Therefore, $D$ is a square, hence $D = t^2$ for some nonnegative integer $t$.    
We have
\begin{align*}
Q(t-b,a) & = \frac{1}{2} \left( a(t-b)^2 + 2b(t-b)a + ca^2 \right)\\
& = \frac{a}{2}(t^2 -b^2 + ac) = 0.
\end{align*}
Recall that $a$ and $c$ are positive integers.  
By Lemma~\ref{Cantor-FP:lemma:QuadraticForm}, 
the quadratic form $Q(x,y)$ is positive-definite 
for $(x,y) \in \N_0^2$, and so $t-b < 0$.  
By Lemma~\ref{Cantor-FP:lemma:abc}, we have $0 \leq t < b \leq 1$ and so $b = 1$.
Therefore, $1-ac = b^2 - ac = t^2 \geq 0$ and so $0 < ac \leq 1$.  
This implies that $a=c=1$, and so 
\[
F(x,y) = \frac{(x+y)^2}{2} +\frac{1}{2}(dx+ey) +f
\]
Moreover, $d\equiv a \equiv 1\pmod{2}$ and  
$e\equiv c \equiv 1\pmod{2}$, that is, $d$ and $e$ are odd integers.

If $d = e$, then $F(x,y) = F(y,x)$ for all $(x,y) \in \N_0^2$ and $F(x,y)$ is not one-to-one.  
Therefore, $d \neq e.$  
If $d > e$, then $d - e = 2g$ for some positive integer $g$, and 
\[
F(x,y) =\frac{(x+y)(x+y+e)}{2}  + gx + f.
\] 
We have $e \neq 0$ because $e$ is odd, and 
\[
F(0,-e) = F(0,0) = f. 
\]
Because $F(x,y)$ is a sorting function, 
we conclude that $e \geq 1$.  Therefore, $F(x,y) \geq f$ for all $(x,y) \in \N_0^2$, 
and so $f=0$.  

If $e \geq 3$, then for all $(x,y) \in \N_0^2\setminus \{ (0,0) \}$ 
we have $x+y \geq 1$ and $F(x,y) \geq (1+e)/2 \geq 2$.  
This means that $F(x,y) \neq 1$ for all $(x,y)\in \N_0^2$, which is absurd.  
Therefore, $e=1$ and
\[
F(x,y) = \frac{(x+y)(x+y+1)}{2} + gx
\] 
for some positive integer $g$.  We have
$F(0,1) = 1$, $F(1,0) = 1+g$, and $F(x,y) \geq 3$ for all $(x,y)\in \N_0^2$ with $x+y \geq 2$.  
If $g\geq 2$, then $F(x,y) \neq 2$ for all  $(x,y)\in \N_0^2$, which is absurd.  
Therefore, $g=1$ and $F(x,y) = C_1(x,y)$ is the first Cantor polynomial.  

Similarly, if $e>d$, then $F(x,y) = C_2(x,y)$ is the second Cantor polynomial.  
This completes the proof.  
\end{proof}

\section{Packing polynomials in sectors}
For every positive real number $\alpha$ we construct the 
\emph{real sector}  \index{real sector}
\[
S(\alpha) = \{ (x,y) \in \R^2 :  0 \leq y \leq \alpha x\}
\]
and the \emph{integer sector}  \index{integer sector}
\[
I(\alpha) = S(\alpha) \cap \N_0^2 
= \{ (x,y) \in \N_0^2 : 0 \leq y \leq \alpha x\}.
\]
The real sector $S(\alpha)$ is the cone with vertex at $(0,0)$ 
generated by the points $(1,0)$ and $(1,\alpha)$.  
A sector is called \emph{rational}  \index{rational sector} if $\alpha$ is a rational number 
and \emph{irrational}  \index{irrational sector} if $\alpha$ is an irrational number.  
Recent work on Cantor polynomials has concentrated on packing polynomials 
in rational sectors.  

\bt[Nathanson~\cite{nath14c}]     \label{CantorSector:theorem:I(r/s)}
Let $r$ and $s$ be relatively prime positive integers 
such that $1 \leq r < s$ and $r$ divides $s-1$.  
Let $d = (s-1)/r$. 
The polynomials
\[
F_{r/s}(x,y) = \frac{r(x-dy)^2}{2} + \frac{ (2-r)x + (dr-2d+2)y }{2}
\]
and
\[
G_{r/s}(x,y) = \frac{r(x-dy)^2}{2} +  \frac{(r+2)x  - (2d+s+1)y}{2}.
\]
are quadratic packing polynomials for the integer sector $I(r/s)$.  
Moreover, for $s \geq 2$, the polynomials
\[
F_{1/s}(x,y) = \frac{(x-(s-1)y)^2}{2} + \frac{x+(3-s)y}{2}
\]
and
\[
G_{1/s}(x,y) = \frac{(x-(s-1)y)^2}{2} + \frac{3x+(1-3s)y}{2}
\]
are the unique quadratic packing polynomials for the integer sector $I(1/s)$.  
\et

Recent work by Stanton~\cite{stan14} and Brandt~\cite{bran14}
has determined explicitly the quadratic packing polynomials for all rational sectors.  

It is an open problem to understand sorting and packing polynomials on irrational sectors.

\def\cprime{$'$} \def\cprime{$'$} \def\cprime{$'$} \def\cprime{$'$}
\providecommand{\bysame}{\leavevmode\hbox to3em{\hrulefill}\thinspace}
\providecommand{\MR}{\relax\ifhmode\unskip\space\fi MR }
\providecommand{\MRhref}[2]{%
  \href{http://www.ams.org/mathscinet-getitem?mr=#1}{#2}
}
\providecommand{\href}[2]{#2}

\end{document}